\documentclass[10pt,colorinlistoftodos]{amsart}

\usepackage{blindtext}
 \usepackage{a4wide}

\usepackage[utf8]{inputenc} 
\usepackage[T1]{fontenc}    
\usepackage[hidelinks]{hyperref}
\usepackage{url}            
\usepackage{booktabs}       
\usepackage{amsfonts}       
\usepackage{amsmath}
\usepackage{amsthm}
\usepackage{amssymb}
\usepackage{nicefrac}       
\usepackage{microtype}      
\usepackage{thmtools}
\usepackage{cleveref}       
\usepackage{lipsum}         
\usepackage{graphicx}
\usepackage{doi}
\usepackage{color}
\usepackage{tikz-cd}
\usepackage{lmodern}
\usepackage{mathtools}
\usepackage[thinc]{esdiff}
\usepackage{cases}
\usepackage{empheq}

\usepackage{cleveref}
\usepackage{lipsum} 
\usepackage{todonotes}
\usepackage{orcidlink}
\usepackage{xcolor}
\usepackage{lineno}

\newcommand{\Q}{\mathbb{Q}}
\newcommand{\R}{\mathbb{R}}
\newcommand{\C}{\mathbb{C}}
\newcommand{\K}{\mathbb{K}}

\newcommand{\F}{\mathbb{F}}

\newtheorem{theorem}{Theorem}[section]
\newtheorem{thm}{Theorem}[section]
\newtheorem{lemma}[thm]{Lemma}
\newtheorem{corollary}[thm]{Corollary}
\newtheorem{prop}[thm]{Proposition}

\newtheorem*{thm*}{Theorem} 

\theoremstyle{definition}

\newtheorem{remark}{Remark}[section]

\title{Galois scalar extension of simple Lie algebras}

\author[A.~Di Bartolo]{Alfonso Di Bartolo\,\orcidlink{0000-0001-5619-2644}}
\address[A.~Di Bartolo]{Dipartimento di Matematica e Informatica,
Università degli Studi di Palermo, Via Archirafi 34, 90123, Palermo, Italy}
\email{alfonso.dibartolo@unipa.it}

\author[F.~P.~Di Fatta]{Francesco Paolo Di Fatta*\,\orcidlink{0009-0008-6197-6126}}
\address[F.~P.~Di Fatta]{Doctoral School on Mathematics and Computational Sciences, 
University of Messina, Catania and Palermo, Italy}
\email{francesco.difatta@studenti.unime.it}
\thanks{*Corresponding author}

\author[G.~La Rosa]{Gianmarco La Rosa\,\orcidlink{0000-0003-1047-5993}}
\address[G.~La Rosa]{Dipartimento di Matematica e Informatica,
Università degli Studi di Palermo, Via Archirafi 34, 90123, Palermo, Italy}
\email{gianmarco.larosa@unipa.it}

\setuptodonotes{color=blue!30}

\definecolor{math}{rgb}{0.0, 0.8, 0.6}


\begin{document}
\maketitle

\begin{abstract}
We investigate the behaviour of simple finite-dimensional Lie algebras under finite Galois field extensions. Given a simple Lie algebra $L$ over a field $\mathbb{F}$ and a finite Galois extension $\mathbb{K}/\mathbb{F}$ of degree $n = [\mathbb{K}:\mathbb{F}]$, we prove that if $L$ is isomorphic to a $\mathbb{K}$-Lie algebra $S$ viewed as an $\mathbb{F}$-algebra, then the scalar extension $L \otimes_{\mathbb{F}} \mathbb{K}$ is isomorphic to the direct sum of $n$ copies of $S$. Furthermore, for cyclic extensions of order a prime number, we show that if no such $\mathbb{K}$-Lie algebra $S$ exists, the extended algebra $L \otimes_{\mathbb{F}} \mathbb{K}$ remains simple. This framework provides a general mechanism to transfer structural and classification results across different ground fields. In addition, we examine the centroids of simple Lie algebras, proving that the dimension of the space of centroids is invariant under finite field extensions. As an application, we explicitly determine the structure of the spaces of centroids for simple Lie algebras over finite fields.
\end{abstract}

\section{Introduction}

The structure and classification of finite-dimensional simple Lie algebras over algebraically closed fields of characteristic zero is a foundational result in modern Lie theory. Over the complex numbers $\mathbb{C}$, simple Lie algebras are fully classified by Dynkin diagrams via their root systems, yielding the four classical families together with the five exceptional Lie algebras \cite{erdmann2006introduction, humphreys2012introduction}. 

The classification of real simple Lie algebras relies on a similar structural mechanism. Any finite-dimensional real Lie algebra $L$ yields a complexification $L \otimes_{\mathbb{R}} \mathbb{C}$, and the complexification of a simple real Lie algebra is semisimple \cite{knapp1996lie}. Consequently, real forms are classified by first studying the corresponding complex Lie algebra and then identifying the non-isomorphic real Lie algebras whose complexifications coincide. This construction illustrates a broader phenomenon: extending the field of scalars often simplifies the internal structure of an algebra, while non-isomorphic algebras over the ground field may become isomorphic after extension.

This behaviour naturally leads to the study of scalar extensions for Lie algebras over arbitrary fields. Let $\mathbb{F}$ be a field and let $\mathbb{K}/\mathbb{F}$ be a finite extension. Given a finite-dimensional Lie algebra $L$ over $\mathbb{F}$, we consider its scalar extension $L \otimes_{\mathbb{F}} \mathbb{K}$.
Even if $L$ is simple over $\mathbb{F}$, the extended algebra is not necessarily semisimple over $\mathbb{K}$. 
An analogous framework is well established in the context of associative algebras. The theory of central simple algebras provides a natural setting for analysing scalar extensions and descent phenomena; see \cite{Pierce} for the general structure theory of associative algebras and \cite{Gille_Szamuely_2006} for the role of Galois cohomology.

However, a fundamental distinction arises when working over ground fields of arbitrary characteristic (see \cite{block1988,Strade1988,Strade2004,Strade2007,Strade2008,Strade2009,Strade2013,Strade2004_2}). In this general setting, standard characteristic-zero tools cannot be applied directly. In particular, although semisimplicity is preserved under the scalar extensions considered here, a semisimple Lie algebra in arbitrary characteristic cannot generally be decomposed into a direct sum of simple ideals. Given a finite Galois extension $\mathbb{K}/\mathbb{F}$ with Galois group $\operatorname{Gal}(\mathbb{K}/\mathbb{F})$, the structure of $L \otimes_{\mathbb{F}} \mathbb{K}$ encodes the descent data back to $L$ via the group.
The scalar extension of Lie algebrasfinds applications in other settings, such as the centroids \cite{jacobson1979lie,Melville1992}, derivations \cite{Burde2021}, biderivations \cite{DiBartolo2026} and also generalised derivations of Lie current algebras \cite{Benkovič2024,Neeb2008,Zusmanovich2005}. Under the standard descent correspondence, any proper ideal of $L \otimes_\F \K$ that remains stable under the action of $G$ corresponds to a proper ideal of $L$. Hence, non-trivial ideals of $L \otimes_\F \K$ are closely linked to the orbit structure induced by the Galois group (\cite{Serre2004,Gille_Szamuely_2006,Knus1998,waterhouse1979}). Understanding the structure of $L \otimes_\F \K$ thus reduces to analysing the interplay between the ideal structure arising after scalar extension and the symmetries enforced by the Galois group.
In this paper, we focus on the case where $\mathbb{K}/\mathbb{F}$ is a cyclic Galois extension of prime degree $p$.

The primary objective of this work is to ccharacterise these structural constraints for scalar extensions of simple Lie algebras. Specifically, we investigate the admissible ideal configurations of $L \otimes_{\mathbb{F}} \mathbb{K}$ when $L$ is simple over $\mathbb{F}$ and $\mathbb{K}/\mathbb{F}$ is cyclic of prime degree. Rather than assuming that the extended algebra decomposes as a direct sum of simple ideals, we leverage the Galois action to describe how its ideals interact and how their configurations are constrained by the simplicity of $L$.
The results presented below establish precise structural limitations in the prime-cyclic case, clarifying the extent to which the field-extension behaviour of a simple Lie algebra over an arbitrary field is constrained by its Galois descent data. The paper is organisedd as follows. In \Cref{sec:preliminaries}, we recall foundational definitions and key background results concerning scalar extensions and Galois action on Lie algebras, including Dedekind's theorem on the linear independence of characters. In \Cref{sec:Main}, we establish our main result: we prove that if $L$ is a simple Lie algebra over $\mathbb{F}$ and $\mathbb{K}/\mathbb{F}$ is a cyclic Galois extension of prime degree, then the scalar extension $L \otimes_{\mathbb{F}} \mathbb{K}$ either remains simple over $\mathbb{K}$ or decomposes as a direct sum of $[\mathbb{K}:\mathbb{F}]$ copies of a simple $\mathbb{K}$-algebra viewed as an $\mathbb{F}$-algebra. Finally, in \Cref{sec:centroids}, we examine the centroid of simple Lie algebras under scalar extensions, providing an a explicit characterisation of these maps in the setting of finite fields.
\section{Preliminaries}\label{sec:preliminaries}

In this section we examine the effect of enlarging the field of scalars for a Lie algebra. Let \(\F\) be a field, \(\K\) be an extension field of \(\F\) and \(L\) a Lie algebra over \(\F\). Since \(\K\), in particular, is a \(\F-\)vector space, we can consider the tensor product \(L \otimes_{\F} \K\). 
If \(c \in \K\), the multiplication by \(c\) is a \(\F-\)linear map from \(\K\) to itself. Thus, for the universal property of the tensor product, we have a well-defined multiplication by \(c\) for \(L \otimes_{\F} \K\) given by \(c (v \otimes a)=v \otimes ca\). Then \(L \otimes_{\F} \K\) is also a \(\K-\)vector space. Now, we want to define a Lie algebra structure over the tensor product: given \(v \otimes a, w \otimes b \in L \otimes_{\F} \K\), we define the Lie bracket as 

\begin{equation*}
    [v \otimes a , w \otimes b]=[v,w] \otimes ab
\end{equation*}

Again, for the universal property of the tensor product, this operation is well-defined, and it arises a Lie algebra structure on \(L \otimes_{\F} \K\) over the field \(\K\).

\begin{prop}\label{prop:uniqueness}
Let \(L\) be a Lie algebra over the field \(\F\) and \(\K\) be a finite extension field of \(\F\), with basis \(B=\big\{e_1,\ldots,e_n\big\}\). Then every elements of \(L \otimes_{\F} \K\) can be written uniquely as \(\sum_{i=1}^n v_i \otimes e_i\), for some \(v_i \in L\).
\end{prop}

\vspace{2mm}

According to Dedekind's Theorem on the linear independence of Galois automorphisms, for any finite Galois extension \(\F \hookrightarrow \K\), the elements of the Galois group are linearly independent as functions over \(\K\).
Formally, we have the following proposition.

\begin{lemma}(Dedekind) \cite{Rotman1998}
Let \(\mathrm{Gal}(\K / \F)=\big\{\tau_1,\ldots,\tau_n\big\}\). Then for any coefficients 
\(\alpha_1, \ldots, \alpha_n\) in \(\K\), the following condition holds:
\begin{equation*}
\sum_{i=1}^n \alpha_i \tau_i = 0 \implies \alpha_i = 0 \quad \forall i \in {1, \dots, n}.
\end{equation*}
\end{lemma}

As a consequence, we have the following result.

\begin{prop} \label{prop:matrix}
Let \(\F \hookrightarrow \K\) be a Galois extension with Galois group \(\mathrm{Gal}(\K / \F)=\big\{\tau_1,\ldots,\tau_n\big\}\) and let \(\big\{a_1,\ldots,a_n\big\}\) be a basis of \(\K\) over \(\F\). Then the following matrix 

\begin{equation*}
\begin{pmatrix}
    \tau_1(a_1) & \tau_2(a_1) & \cdots &\tau_n(a_1)\\
    \tau_1(a_2) & \tau_2(a_2) & \cdots & \tau_n(a_2)\\
    \vdots & \vdots & \ddots & \vdots \\
    \tau_1(a_n) & \tau_2(a_n) & \cdots &\tau_n(a_n)
\end{pmatrix}
\end{equation*}
is invertible over \(\K\).
\end{prop}

\begin{proof}
Assume, for the sake of contradiction, that the determinant of the matrix is zero. This implies that its columns are linearly dependent over 
\(\K\); thus, there exist coefficients 
\(\lambda_1,\ldots,\lambda_n\in \K\), not all zero, such that 

\begin{equation*}
    \sum_{k=1}^n \lambda_k \tau_k(a_j)=0
\end{equation*}
for all \(j=1,\ldots, n\). 
Let \(x=\sum_{l=1}^n \beta_l a_l \in \K\)
be an arbitrary element of the field expressed in the 
\(\F\)-basis \(\big\{a_l\big\}\), where \(\beta_l \in \F\). By the \(\F\)-linearity of the automorphisms, we obtain the following identity:

\begin{equation*}
    \sum_{k=1}^n \lambda_k \tau_k(x)=\sum_{k=1}^n \lambda_k \tau_k(\sum_{l=1}^n \beta_l a_l)=
    \sum_{l=1}^n \beta_l\sum_{k=1}^n \lambda_k \tau_k(a_l)=\sum_{l=1}^n \beta_l\cdot 0=0.
\end{equation*}

Since this holds for every \(x \in \K\), the linear combination of automorphisms vanishes identically as a function. This provides a direct contradiction to Dedekind's Theorem on the linear independence of Galois automorphisms. Consequently, the columns must be linearly independent, and the matrix is necessarily non-singular.
\end{proof}

\section{Main results}\label{sec:Main}
Let \(\mathfrak{s}\) be a Lie algbera over the field \(\K\), for any \(\sigma \in \mathrm{Gal}(\K / \F)\) we define  the conjugate Lie algebra \(\mathfrak{s}^\sigma\). We fix a basis \(\big\{s_1,\ldots,s_m\big\}\) over \(\K\). The sum of vector on \(\mathfrak{s}^\sigma\) is the same as \(\mathfrak{s}\), the scalar multiplication \(\cdot_\sigma\) on \(\mathfrak{s}^\sigma\) is defined as 
\begin{equation*}
    \lambda \cdot_\sigma v= \sigma(\lambda) v
\end{equation*}
and the constant structures of the bracket \([\cdot,\cdot]_{\mathfrak{s}^\sigma}\) on \(\mathfrak{s}^\sigma\) are \(\big\{\sigma(c_{ij}^k)\big\}_{ijk}\)
where \(\big\{c_{ij}^k\big\}_{ijk}\) are the constant structures of \(\mathfrak{s}\) (i.e., \([s_i,s_j]=\sum_{i=1}^m c_{ij}^k s_k\)). The bracket \([\cdot,\cdot]_{\mathfrak{s}^\sigma}\) is indeed a Lie bracket. In fact by \(c_{ii}^k=0\) follows \(\sigma(c_{ii}^k)=\sigma(0)=0\) and by \(\sum_{l=1}^m \left( c_{ij}^{l} c_{lk}^{n} + c_{jk}^{l} c_{li}^{n} + c_{ki}^{l} c_{lj}^{n} \right)=0\) it follows
\begin{equation*}
    0=\sigma(\sum_{m} \left( c_{ij}^{\;\;m} c_{mk}^{\;\;n} + c_{jk}^{\;\;m} c_{mi}^{\;\;n} + c_{ki}^{\;\;m} c_{mj}^{\;\;n} \right)) = \sum_{m} \left( \sigma(c_{ij}^{\;\;m}) \sigma(c_{mk}^{\;\;n}) + \sigma(c_{jk}^{\;\;m}) \sigma(c_{mi}^{\;\;n}) + \sigma(c_{ki}^{\;\;m}) \sigma(c_{mj}^{\;\;n}) \right).
\end{equation*}
We define \(\Tilde{\sigma}: s \to s^\sigma\) such that
\begin{equation*}
    \Tilde{\sigma}(\sum_{i=1}^m\lambda_i s_i)=\sum_{i=1}^m \sigma(\lambda_i) s_i=\sum_{i=1}^m \lambda_i \cdot_\sigma s_i
\end{equation*}
It is obvious a homomorphism of Lie algebras, it is \(\K\)-linear and preserve the bracket. Furthermore, by the bijective of \(\sigma\), it follows that \(\Tilde{\sigma}\) is isomorphism of Lie algebras.

\begin{remark}
We observe that the map \(\tilde{\sigma}\) depend on the choice of the basis, if we take another basis \(\big\{v_1,\ldots,v_m\big\}\) with \([v_i,v_j]=\sum_{k=1}^n f_{ij}^k v_k\), it is not true that \([v_i,v_j]_{\mathfrak{s}^\sigma}=\sum_{k=1}^n \sigma(f_{ij}^k) v_k\).
\end{remark}
Let \(L\) be a Lie algebra over \(\F\) with basis $\{v_1, \dots, v_m\}$ of $L$. By \Cref{prop:uniqueness} it follows that $\{v_i \otimes 1\}_{i=1, \dots, m}$ is a $\K$-basis of $L \otimes_\F \K$. Then we have
\(\Tilde{\sigma}: L \otimes_\F \K \to (L \otimes_\F \K)^\sigma\) such that

\begin{equation*}
    \Tilde{\sigma}(\sum_{i=1}^m\lambda_i \left(v_i \otimes 1\right))=\sum_{i=1}^n \sigma(\lambda_i) (v_i \otimes 1)=\sum_{i=1}^n v_i \otimes \sigma(\lambda_i) 
\end{equation*}

\begin{prop}\label{prop:ideal for sigma}
    Let \(\mathfrak{s}\) be a Lie algebra over \(\K\), let \(\big\{s_1,\ldots,s_m\big\}\) be a basis over \(\K\) such that the costant structures are in \(\F\), i.e. \([s_i,s_j]=\sum_{k=1}^m c_{ij}^k s_k\) with \(c_{ij}^k \in \F\) for all \(i,j,k\). Then \(\tilde{\sigma}([x,y])=[\tilde{\sigma}(x),\tilde{\sigma}(y)]\) for all \(x,y \in \mathfrak{s}\).
\end{prop}

\begin{proof}
Let \(x=\sum_{i=1}^m \alpha_i s_i\) and \(y=\sum_{j=1}^m \beta_j v_j\). Then we have
\begin{align*}
    \tilde{\sigma}([x,y])&=\tilde{\sigma}(\sum_{i,j=1}^m \alpha_i \beta_j [s_i,s_j])\\
    &=\tilde{\sigma}(\sum_{k=1}^m(\sum_{i,j=1}^m \alpha_i \beta_j c_{ij}^k)s_k)\\
    &=\sum_{k=1}^m \sigma(\sum_{i,j=1}^m \alpha_i \beta_j c_{ij}^k) s_k\\
    &=\sum_{k=1}^m (\sum_{i,j=1}^m(\sigma( \alpha_i)\sigma (\beta_j) \sigma(c_{ij}^k)) s_k\\
    &=\sum_{k=1}^m (\sum_{i,j=1}^m(\sigma( \alpha_i)\sigma (\beta_j) c_{ij}^k) s_k
\end{align*}
because \(c_{ij}^k\) are elements of \(\F\). But we also have
\begin{align*}
    [\tilde{\sigma}(\sum_{i=1}^m \alpha_i s_i),\tilde{\sigma}(\sum_{j=1}^m \beta_j s_j)]&=[\sum_{i=1}^m \sigma(\alpha_i) s_i,\sum_{j=1}^m \sigma(\beta_j) s_j]\\
    &=\sum_{i,j=1}^m \sigma(\alpha_i)\sigma(\beta_j)[s_i,s_j]\\
    &=\sum_{k=1}^m (\sum_{i,j=1}^m(\sigma( \alpha_i)\sigma (\beta_j) c_{ij}^k) s_k
\end{align*}
\end{proof}

\begin{remark}
    An example is the algebra \(L \otimes_\F \K\), indeed if \(\big\{v_1,\ldots,v_n\big\}\) is a \(\F-\) basis of \(L\), then \(\big\{v_i \otimes 1\big\}_{i=1,\ldots,n}\) is a \(\K\)- basis of \(L \otimes_\F \K \) and \([v_i \otimes 1, v_j \otimes 1]=[v_i,v_j] \otimes 1 = (\sum_{k=1}^n c_{ij}^k v_k) \otimes 1= \sum_{k=1}^n c_{ij}^k (v_k \otimes 1)\) with \(c_{ij}^k \in \F\).
\end{remark}
\begin{corollary}
    Let \(L\) be a Lie algebra over \(\F\), let \(\K\) be a Galois extension of \(\F\) and let \(I\) be an ideal of \(L \otimes_\F \K\). then \(I_\sigma=\tilde{\sigma}(I)\) is also an ideal of \(L \otimes_\F \K\).
\end{corollary}
\begin{proof}
   \(\tilde{\sigma}: L \otimes_\F \K \to (L \otimes_\F \K)^\sigma \) is an isomorphism of Lie algebra, so \(I_\sigma\) is an ideal of \((L \otimes_\F \K)^\sigma\).
    \(I_\sigma\) is closed under addition because \(I\) is an ideal and the two addition rules coincide. it is closed under scalar multiplication. Indeed let \(\beta \in K\) and \(x \in I\). then there exists \(\alpha \in \K\) such that \(\beta=\sigma(\alpha)\). Hence
    \begin{equation*}
        \beta x= \sigma(\alpha) x=\alpha \cdot_\sigma x \in I_\sigma
    \end{equation*}
    So \(I_\sigma\) is a subspace of \(L \otimes_\F \K\). Lastly, let \(x \in L \otimes_\F \K\) and \(y \in I_\sigma\). Then there exists \(x' \in L \otimes_\F \K\) and \(y' \in I\) such that \(x=\tilde{\sigma}(x')\) and \(y=\tilde{\sigma}(y')\). Hence according to the previous remark \([x,y]=[\tilde{\sigma}(x'),\tilde{\sigma}(y')]=\tilde{\sigma}([x',y'])\in \tilde{\sigma}(I)=I_\sigma\).
\end{proof}
\begin{lemma}\label{lem:independence}

Let $\F \hookrightarrow \K$ be a finite Galois extension of degree $[\K : \F] = n$, with $\mathrm{Gal}(\K/\F) = \{\sigma_1, \dots, \sigma_n\}$. Let $\{a_1, \dots, a_n\}$ be a basis of $\K$ over $\F$, and let $\mathfrak{s}$ be a vector space over $\K$. Given $y_1, \dots, y_n \in \mathfrak{s}$, if 
\[ \sum_{i=1}^n \sigma_l(a_i) y_i = 0 \]
for every $l = 1, \dots, n$, then $y_1 = \dots = y_n = 0$.
\end{lemma}

\begin{proof}
 Let \(\big\{s_1,\ldots,s_m\big\}\) be a basis of \(\mathfrak{s}\) over \(\K\). Then there exist 
 \(\beta_{ij}\in \K\) such that \(y_i=\sum_{j=1}^m \beta_{ij}s_j\). Then the identity became
\begin{equation*}
     0= \sum_{i=1}^n \sigma_l(a_i) y_i = \sum_{j=1}^m\left(\sum_{i=1}^n \sigma_l(a_i) \beta_{ij}\right)s_j 
\end{equation*}
for every \(l=1,\ldots,n\). Then we have
\begin{equation*}
     \sum_{i=1}^n \sigma_l(a_i) \beta_{ij}=0.
\end{equation*}
By \Cref{prop:matrix} the following system
\[
\begin{pmatrix}
\sigma_1(a_1) & \sigma_1(a_2) & \cdots & \sigma_1(a_n) \\
\sigma_2(a_1) & \sigma_2(a_2) & \cdots & \sigma_2(a_n)\\
\vdots & \vdots & \ddots & \vdots \\
\sigma_n(a_1) & \sigma_n(a_2) & \cdots & \sigma_n(a_n) \\
\end{pmatrix}
\begin{pmatrix}
\beta_{1j} \\
\beta_{2j} \\
\vdots\\
\beta_{nj}
\end{pmatrix}
=
\begin{pmatrix}
0 \\
0 \\
\vdots\\
0
\end{pmatrix}
\]
has only the solutions \(\beta_{1j}=\beta_{2j}=\ldots=\beta_{nj}=0\) for every \(j\).
\end{proof}

\begin{lemma}\label{lem:invarianza}
Let $L$ be a simple Lie algebra over $\F$ and let $\F \hookrightarrow \K$ be a finite Galois extension. An ideal \(I\) of \(L \otimes_\F \K\) satisfy \(I_\sigma=I\) for every \(\sigma \in \mathrm{Gal}(\K /\F)\) if and only if \(I=0_{L \otimes_\F \K}\) or \(I=L \otimes_\F \K\).
\end{lemma}

\begin{proof}
    Let \(\big\{a_1,\ldots,a_n\big\}\) be a basis of \(\K\) over \(\F\) and let \(\mathrm{Gal}(\K / \F)=\big\{\sigma_1,\ldots,\sigma_n\big\}\). If \(I\) is an improper ideal, the statement is clear. Let \(I\) be a non zero ideal of \(L \otimes_\F \K\) and let \(\sum_{i=1}^p x_i \otimes a_i \in I\) be a nonzero vector. So there exists \(x_l \neq 0\) for some $k$. Since \(I\) is invariant for every \(\tilde{\sigma}_j\), \(\sum_{i=1}^p x_i \otimes \sigma_j(a_i) \in I\). By \Cref{prop:matrix} the following linear system

\[
\begin{pmatrix}
\sigma_1(a_1) & \sigma_2(a_1) & \cdots & \sigma_l(a_1) & \cdots & \sigma_n(a_1) \\
\sigma_1(a_2) & \sigma_2(a_2) & \cdots & \sigma_l(a_2) & \cdots  & \sigma_n(a_2)\\
\vdots & \vdots & \ddots & \vdots & \ddots  & \vdots\\
\sigma_1(a_l) & \sigma_2(a_l) & \cdots & \sigma_l(a_l) & \cdots & \sigma_n(a_l) \\
\vdots & \vdots & \ddots & \vdots & \ddots  & \vdots\\
\sigma_1(a_n) & \sigma_2(a_n) & \cdots & \sigma_l(a_n) & \cdots & \sigma_n(a_n) \\
\end{pmatrix}
\begin{pmatrix}
k_1 \\
k_2 \\
\vdots\\
k_l\\
\vdots\\
k_n
\end{pmatrix}
=
\begin{pmatrix}
0 \\
0 \\
\vdots\\
a_l\\
\vdots\\
0
\end{pmatrix}
\]
has a unique solution \(k_1,\ldots,k_n \in \K\). So we have
\begin{equation*}
    x_l \otimes a_l= \sum_{i=1} x_i \otimes \left(\sum_{r=1}^n k_r \sigma_r(a_i)\right)=\sum_{r=1}^n k_r\left(\sum_{i=1} x_i \otimes \sigma_r(a_i)\right)\in I
\end{equation*}

The space generated by the elements of the form \([z_1,[z_2,[\ldots,[z_n,x_l]]\ldots]]\) gives an ideal of \(L\). Since \(L\) is simple and \(x_l \neq 0\), this ideal coincides with \(L\). Therefore, if we apply \(z_i \otimes 1\) to \(x_l \otimes a_l\) we obtain that \(L \otimes a_l \subseteq I\). Hence \(\frac{a_i}{a_l}(L \otimes a_l)=L \otimes a_i \in L\) for every \(i\) and therefore \(I=L \otimes_\F \K\).
\end{proof}

In \cite{knapp1996lie} the authors give a caractherization of the algebras \(L \otimes_\R \C\).

\begin{theorem}[p.~407, Theorem~6.94, \cite{knapp1996lie}]\label{complexification_real}
    Let \(L_0\) be a simple Lie algebra over $\R$ and let $L$ be its complexification. Then there are just two possibilities:

    \begin{itemize}
        \item $L_0$ is a complex, i.e., is of the form $\mathfrak{s}_\R$ for some complex $\mathfrak{s}$, and then $L$ is $\C$ isomorphic to $\mathfrak{s}\oplus \mathfrak{s}$,
        \item $L_0$ is not complex, and then $L$ is simple over $\C$.
    \end{itemize}
\end{theorem}

\begin{remark}
    The previous result is equivalent to stating that, once a real simple Lie algebra \(L\) is fixed, we have.
        \begin{itemize}
        \item If there exists a Lie algebra isomorphism \(\varphi: L \to \mathfrak{s}_\R\)  over \(\R\), where $\mathfrak{s}$ is a complex Lie algebra, then \(L \otimes_{\R} \C \cong \mathfrak{s} \oplus \mathfrak{s}\) is a semisimple Lie algebra and it is the direct sum of two simple Lie algebras.
        \item Otherwise, 
        \(L \otimes_{\R} \C \) is simple.
    \end{itemize} 
\end{remark}

From now on, we consider only cyclic Galois extension \(K\) of \(\F\) with \([\K:\F]=p\) where \(p\) is a prime number. So the Galois group
\(\mathrm{Gal}(\K / \F)\) has the form
\begin{equation*}
    \mathrm{Gal}(\K / \F)=\big\{\sigma, \sigma^2, \ldots, \sigma^{p-1},\sigma^p=\mathrm{id}\big\}.
\end{equation*}
and every \(\sigma^i\) generates the entire group for all \(i=1,\ldots,p-1\). 

We now generalise the result in \cite{knapp1996lie} to \(L \otimes_\F \K\) where \(\F \hookrightarrow \K\) is a cyclic Galois extension.

\begin{theorem}\label{th:isomorphism}
Let $L$ be a simple Lie algebra over $\F$, and let $\F \hookrightarrow \K$ be a cyclic 
Galois extension of order $p$. Let $L \otimes_\F \K$ denote its scalar 
extension. Then, if there exists an $\F$-isomorphism of Lie algebras $\varphi: L \to \mathfrak{s}_\F$, 
where $\mathfrak{s}_\F$ denotes the \(\K-\)Lie algebra $\mathfrak{s}$ regarded as an 
$\F$-algebra via restriction of scalars, then the tensor product $L \otimes_\F \K$ 
is a semisimple Lie algebra isomorphic to the direct sum of $p$ copies of $\mathfrak{s}$:
\begin{equation*}
    L \otimes_\F \K \cong \underbrace{\mathfrak{s} \oplus \dots \oplus \mathfrak{s}}_{p \text{ times}},
\end{equation*}
\end{theorem}

\begin{proof}
Let $\mathrm{Gal}(K/F) = \{\sigma, \sigma^2, \dots, \sigma^{p-1},\sigma^p=\mathrm{id}\}$ be the Galois group, let \(\big\{a_1,\ldots,a_p\big\}\) be a basis of \(\K\) over \(\F\) and let $\varphi: L \to \mathfrak{s}_\F$ be an $\F$-isomorphism of Lie algebras. Given an $\F$-basis $\{v_1, \dots, v_m\}$ of $L$, by \Cref{prop:uniqueness} it follows that $\{v_i \otimes 1\}_{i=1, \dots, m}$ is a $\K$-basis of $L \otimes_\F \K$. We define \(F : L \otimes_\F \K \to  \mathfrak{s}^\sigma \oplus \mathfrak{s}^{\sigma^2} \oplus \cdots \oplus \mathfrak{s}^{\sigma^p}\) with
\begin{equation*}
           F( v_i \otimes 1)=(\varphi(v_i),\varphi(v_i),\ldots,\varphi(v_i))
\end{equation*}
       and we extend by \(\K-\)linearity. By \(\F-\)linearity of \(\varphi\) we obtain
       \begin{equation*}
           F(x \otimes 1)=(\varphi(x),\varphi(x),\ldots,\varphi(x))
       \end{equation*}
       for every \(x \in L\). Then 
       \begin{align*}
           F(\sum_{i=1}^p x_i \otimes a_i)&=\sum_{i=1}^p a_l F(x_i \otimes 1)=\sum_{i=1}^n a_i (\varphi(x_i),\varphi(x_i),\ldots,\varphi(x_i))\\
           &=\sum_{i=1}^p a_i (\varphi(x_i),\varphi(x_i),\ldots,\varphi(x_i))\\
           &=\sum_{i=1}^p (a_i \cdot_\sigma \varphi(x_i), a_i \cdot_{\sigma^2} \varphi(x_i),\ldots, a_i \cdot_{\sigma^p} \varphi(x_i))\\
           &=\sum_{i=1}^p (\sigma(a_i) \varphi(x_i), \sigma^2(a_i) \varphi(x_i),\ldots, \sigma^p(a_i) \varphi(x_i))\\
           &=(\sum_{i=1}^p \sigma(a_i) \varphi(x_i), \sum_{i=1}^n \sigma^2(a_i) \varphi(x_i),\ldots, \sum_{i=1}^p \sigma^p(a_i) \varphi(x_i))
       \end{align*}
      Since \(\varphi\) is a homomorphism of Lie algebras, \(F\) is likewise a 
      homomorphism. To establish this assertion, it is sufficient to check the 
      condition solely on the basis vectors \(\big\{v_i \otimes 1\big\}\).
       \begin{align*}
           F([v_i \otimes 1, v_j \otimes 1])&=F([v_i,v_j] \otimes 1)\\
           &=(\varphi([v_i,v_j]),\varphi([v_i,v_j]),\ldots,\varphi([v_i,v_j]))\\
           &=([\varphi(v_i),\varphi(v_j)],[\varphi(v_i),\varphi(v_j)],\ldots,[\varphi(v_i),\varphi(v_j)])\\
           &=[(\varphi(v_i),\varphi(v_i),\ldots,\varphi(v_i)),(\varphi(v_j),\varphi(v_j),\ldots,\varphi(v_j))]\\
           &=[F(v_i \otimes 1),F(v_j \otimes 1)].
       \end{align*}Let \(F(\sum_{i=1}^p x_i \otimes a_i)=0\). Then \(\sum_{i=1}^p \sigma^j(a_i)\varphi(x_i)=0\). According to \Cref{lem:independence} \(\varphi(x_i)=0\) for every \(i\). Since \(\varphi\) is an isomorphism, \(x_i=0\) for all \(i\). Hence \(F\) is injective. We also have
       \begin{align*}
           \dim_\K L \otimes_\F \K &= \dim_\F L \\
           &=\dim_\F \mathfrak{s} = [\K : F] \dim_\K \mathfrak{s}\\
           &= \dim_\K(\mathfrak{s}^\sigma \oplus \mathfrak{s}^{\sigma^2} \oplus \cdots \oplus \mathfrak{s}^{\sigma^n}).
       \end{align*}
        So \(F\) is an isomorphism. Lastly, \(\mathfrak{s}^{\sigma^j} \cong \mathfrak{s}\) via the map \(\tilde{\sigma^j}\), so
        \begin{equation*}
            L \otimes_\F \K \cong \underbrace{\mathfrak{s} \oplus \dots \oplus \mathfrak{s}}_{p \text{ times}},
        \end{equation*} 
\end{proof}

\begin{theorem}\label{th:simple}
Let $L$ be a simple Lie algebra over $\F$, and let $\F \hookrightarrow \K$ be a cyclic 
Galois extension of order $p$. Let $L \otimes_\F \K$ denote its scalar 
extension. Then if there no exists an $\F$-isomorphism of Lie algebras $\varphi: L \to \mathfrak{s}_\F$, 
\(L \otimes_\F \K\) is simple.
\end{theorem}

\begin{proof}
Let $I$ be a non-zero ideal of $L \otimes_\F \K$. The sum of ideals $I_\sigma + I_{\sigma^2} + \dots + I_{\sigma^p}$ is invariant under all automorphisms $\tilde{\sigma}^j$. Thus, according to \Cref{lem:invarianza}, it must be either the zero ideal or the entire space. Since $I \neq 0$, it follows that $L \otimes_\F \K = I_\sigma + I_{\sigma^2} + \dots + I_{\sigma^p}$. Suppose that exists \(i \neq j\) such that \(I_{\sigma^i}=I_{\sigma^j}\). Then we have
$I = \tilde{\sigma}^{p-i}(I_{\sigma^i}) = \tilde{\sigma}^{p-i}(I_{\sigma^j})=I_{\sigma^{j-i}}$. By setting $k = j-i \neq 0$, we obtain:
        \begin{equation*}
            I = I_{\sigma^k} = \tilde{\sigma}^k(I) = \tilde{\sigma}^k(I_{\sigma^k}) = I_{\sigma^{2k}}.
        \end{equation*}
    By iterating this process, it follows that
        \begin{equation*}
            I = I_{\sigma^{rk}}
        \end{equation*}
    for every $r \in \mathbb{Z}$. Since $\sigma^k$ generates the Galois group, we conclude that $I = I_{\sigma^j}$ for every $j=1, \dots, p$. Consequently, $L \otimes_\F \K = I_\sigma + I_{\sigma^2} + \dots + I_{\sigma^p} = I$.
So we can suppose that 
\(I_{\sigma^i}\neq I_{\sigma^j}\) for all \(i \neq j\). If
$I_{\sigma^i} \cap I_{\sigma^j} = \{0\}$ for all $i \neq j$, then the sum is direct, i.e., $L \otimes_\F \K = I_\sigma \oplus I_{\sigma^2} \oplus \dots \oplus I_{\sigma^p}$. Let $\pi_j: L \otimes_\F \K \to I_{\sigma^j}$ be the natural projection.
By composing the canonical injection with the projection on \(I_{\sigma_j}\), we obtain the $\mathbb{F}$-homomorphism
\begin{equation*}
    L \xhookrightarrow{\iota} L \otimes_{\mathbb{F}} K \xrightarrow{\pi_j} I_{\sigma^j}.
\end{equation*}
The composition $\pi_j \circ \iota$ must have a trivial kernel for at least one $j$. Since $L$ is a simple algebra, there exists $l \in \{1, \dots, n\}$ such that $\ker(\pi_l \circ \iota) = 0$. Moreover \(\dim_\F L=\dim_\K L \otimes_\F \K=p \dim_\K I_{\sigma^j}\) (we have effectively \(p\) distinct copies that does not intersect) \(=\dim_\F I_{\sigma^j}\) (the dimensions of \(I_{\sigma^i}\) is equal to the dimension of \(I\) for every \(i\)).
Therefore, $\pi_j \circ \iota$ is an $\mathbb{F}$-isomorphism.
This leads to a contradiction, as $L$ cannot be $\F$-isomorphic to a Lie algebra defined over $\K$ by hypothesis. So there exists a nonzero intersection \(I_{\sigma^i} \cap I_{\sigma^j}\) with \(i<j\). We have two possibilities:
\begin{itemize}
    \item If $\dim_K (I_{\sigma^i} \cap I_{\sigma^j}) = \dim_K I_{\sigma^i} = \dim_K I_{\sigma^j}$, it follows that $I_{\sigma^i} = I_{\sigma^j}$. This cannot be because we suppose \(I_{\sigma^i} \neq I_{\sigma^j}\) for every \(i \neq j\).
    \item Suppose now that $\dim_K (I_{\sigma^i} \cap I_{\sigma^j}) < \dim_K I_{\sigma^i}$. By setting $J = I_{\sigma^i} \cap I_{\sigma^j}$, we see that $J$ is an non zero ideal and the sum $J_\sigma + J_{\sigma^2} + \dots + J_{\sigma^p}$ is an ideal invariant under all automorphisms $\tilde{\sigma}^j$. Thus, according to \Cref{lem:invarianza}, it must be either the zero ideal or the entire space. Since $J \neq \{0\}$, then $L \otimes_\F \K = J_\sigma + J_{\sigma^2} + \dots + J_{\sigma^p}$. We can then apply to $J$ the same reasoning previously applied to $I$. We can suppose that \(J_{\sigma^i}\neq J_{\sigma^j}\) otherwise \(L \otimes_\F \K =J \subseteq I\). If every intersection is zero we have again an \(\F-\)isomorphism with a \(K\)-algebra , if there exists a non-zero intersection $J_{\sigma^i} \cap J_{\sigma^j}$ that is strictly contained in $J_{\sigma^i}$, we set $H = J_{\sigma^i} \cap J_{\sigma^j}$ and repeat the previous argument. Since the Lie algebra is finite-dimensional, this descent process must eventually terminate at the zero ideal, yielding a contradiction with the fact that every intersection step is non-trivial.
    Consequently, one of these terms must coincide with the entire space. Since this space is contained in $I_{\sigma^i}$, we can conclude that $L \otimes_\F \K = I_{\sigma^i}$ and therefore $L \otimes_\F \K = I$.
\end{itemize}
\end{proof}

We observe that \Cref{th:isomorphism} holds for an arbitrary Galois extension. Conversely, \Cref{th:simple} does not extend to a general Galois extension, as evidenced by the following remark.

\begin{remark}
Let $L$ be a Lie algebra over a field $\mathbb{F}$, and let $\mathbb{K}/\mathbb{F}$ be a field extension. If the scalar extension $L \otimes_{\mathbb{F}} \mathbb{K}$ is simple, then $L$ is necessarily simple. Indeed, suppose toward a contradiction that there exists a proper ideal $I \subsetneq L$. Then $I \otimes_{\mathbb{F}} \mathbb{K}$ would constitute a proper ideal of $L \otimes_{\mathbb{F}} \mathbb{K}$, yielding a contradiction.

Now we consider the indefinite orthogonal Lie algebra $\mathfrak{so}(1, 3; \mathbb{\F})$ over the field $\mathbb{\F}$, i.e.
\begin{equation*}
    \mathfrak{so}(1,3; \F)=
    \left\{\begin{pmatrix}
    0 & a & b & c \\
    a & 0 & -d & -e \\
    b & d & 0 & -f \\
    c & e & f & 0 \\
    \end{pmatrix}\mid a,b,c,d,e,f \in \F\right\}.
\end{equation*}

In particular, if \(\F=\Q\) then \(\mathfrak{so}(1, 3; \mathbb{Q})\) is simple because its realification, $\mathfrak{so}(1, 3; \mathbb{Q}) \otimes_{\mathbb{Q}} \mathbb{R} \cong \mathfrak{so}(1, 3)$, is the Lorentz Lie algebra, which is known to be simple.

Consider the Galois extension $\mathbb{Q}(\sqrt{2}, i)$ over $\mathbb{Q}$. Since $\dim_{\mathbb{Q}}(\mathfrak{so}(1, 3; \mathbb{Q})) = 6$, and any Lie algebra over $\mathbb{Q}(\sqrt{2}, i)$ must have a $\mathbb{Q}$-dimension that is a multiple of $[\mathbb{Q}(\sqrt{2}, i) : \mathbb{Q}] = 4$, an isomorphism of $\mathbb{Q}$-algebras is impossible. However, the scalar extension $\mathfrak{so}(1, 3; \mathbb{Q}) \otimes_{\mathbb{Q}} \mathbb{Q}(\sqrt{2}, i) \cong \mathfrak{so}(1,3; \mathbb{Q}(\sqrt{2}, i))$ is no longer simple and decomposes into two isomorphic direct summands.

 To explicitly show the decomposition, let $J_l$ and $K_l$ (for $l=1,2,3$) be the generators of rotations and boosts, respectively, for $\mathfrak{so}(1,3; \mathbb{Q})$ i.e.
\begin{equation*}
    J_1=\begin{pmatrix}
        0 & 0 & 0 & 0\\
        0 & 0 & 0 & 0\\
        0 & 0 & 0 & -1\\
        0 & 0 & 1 & 0\\
    \end{pmatrix}
    \qquad
    J_2=\begin{pmatrix}
        0 & 0 & 0 & 0\\
        0 & 0 & 0 & 1\\
        0 & 0 & 0 & 0\\
        0 & -1 & 0 & 0\\
    \end{pmatrix}
    \qquad
    J_3=\begin{pmatrix}
        0 & 0 & 0 & 0\\
        0 & 0 & -1 & 0\\
        0 & 1 & 0 & 0\\
        0 & 0 & 0 & 0\\
    \end{pmatrix}
\end{equation*} 
\begin{equation*}
    K_1=\begin{pmatrix}
        0 & -1 & 0 & 0\\
        -1 & 0 & 0 & 0\\
        0 & 0 & 0 & 0\\
        0 & 0 & 0 & 0\\
    \end{pmatrix}
    \qquad
        K_2=\begin{pmatrix}
        0 & 0 & -1 & 0\\
        0 & 0 & 0 & 0\\
        -1 & 0 & 0 & 0\\
        0 & 0 & 0 & 0\\
    \end{pmatrix}
    \qquad
        K_3=\begin{pmatrix}
        0 & 0 & 0 & -1\\
        0 & 0 & 0 & 0\\
        0 & 0 & 0 & 0\\
        -1 & 0 & 0 & 0\\
    \end{pmatrix}.
\end{equation*}
The brackets are 
\begin{equation*}
    [J_i,J_j]=\epsilon_{ijk}J_k 
    \qquad
    [J_i,K_j]=\epsilon_{ijk}K_k
    \qquad
    [K_i,K_j]=-\epsilon_{ijk}J_k
\end{equation*}
 where \(e_{ijk}\) is the antisymmetric Levi-Civita tensor defined as \(\epsilon_{123}=1\) and the rule that the sign flips when you swap any two indices.  
 Over the extension field $\mathbb{K} = \mathbb{Q}(\sqrt{2}, i)$, on \(\mathfrak{so}(1,3;\Q(\sqrt{2},i))\) we can define the complexified generators:
\[
A_l = \frac{1}{2}(J_l + i K_l), \quad B_l = \frac{1}{2}(J_l - i K_l)
\]
A direct computation of the commutation relations shows that $[A_l, B_m] = 0$ for all $l,m$. For more details, see \cite[p.~106]{Schwartz_2013}.

\end{remark}

\begin{remark}\label{solvability}
  
Now, we suppose that $L$ is a simple Lie algebra over $\mathbb{F}$ and let 
$\mathbb{F} \hookrightarrow \mathbb{K}$ be a solvable Galois extension; that is, 
the group $G := \mathrm{Gal}(\mathbb{K}/\mathbb{F})$ is solvable. Then there 
exists a subnormal series
\begin{equation*}
    \{\mathrm{id}\} = G_t \lhd G_{t-1} \lhd \dots \lhd G_0 = G
\end{equation*}
such that each quotient $G_{i-1}/G_i$ is a cyclic group of prime order. 
By the Fundamental Theorem of Galois Theory, there exists a tower of finite extensions
\begin{equation*}
    \mathbb{F} = \mathbb{F}_0 \hookrightarrow \mathbb{F}_1 \hookrightarrow \dots \hookrightarrow \mathbb{F}_t = \mathbb{K}
\end{equation*}
such that each $\mathbb{F}_i \hookrightarrow \mathbb{F}_{i+1}$ is a Galois extension of prime degree 
\(p_i\).
Then 
\begin{equation*}
    L \otimes_\F \K = L \otimes_{\F_0} \F_t \cong (L \otimes_{\F_0} \F_1) \otimes_{\F_1} \F_t.
\end{equation*}
Therefore
\begin{itemize}
    \item If there no exists a \(\F_0\)-isomorphism \(\varphi: L \to \mathfrak{s}_{F_0}\), where \(s\) is a \(\F_1\)-Lie algebra, then \(\mathfrak{s}=L \otimes_{\F_0} \F_1\) and we can consider
    \begin{equation*}
        \mathfrak{s} \otimes_{\F_1} \F_t=(\mathfrak{s} \otimes_{\F_1} \F_2)\otimes_{\F_2} \F_t
    \end{equation*}
    So we need to analyze \(\mathfrak{s} \otimes_{\F_1} \F_2\) using \Cref{th:isomorphism,th:simple}.
    \item if there exists a \(\F_0\)-isomorphism \(\varphi: L \to \mathfrak{s}_{F_0}\), where \(s\) is a \(\F_1\)-Lie algebra, then
    \begin{equation*}
            L \otimes_{\F_0} \F_1 \cong \underbrace{\mathfrak{s} \oplus \dots \oplus \mathfrak{s}}_{p_0 \text{ times}}.
    \end{equation*} 
    Hence 
    \begin{equation*}
        L \otimes_\F \K \cong \underbrace{\mathfrak{s} \otimes_{\F_1} \F_t \oplus \dots \oplus \mathfrak{s}\otimes_{\F_1} \F_t}_{p_0 \text{ times}}.
    \end{equation*}
    We have that \(\mathfrak{s}\) is a simple Lie algebra over \(\F_1\). For every piece can now consider
    \begin{equation*}
        \mathfrak{s} \otimes_{\F_1} \F_t=(\mathfrak{s} \otimes_{\F_1} \F_2)\otimes_{\F_2} \F_t
    \end{equation*}
    and we reapply \Cref{th:isomorphism,th:simple} to every piece \(\mathfrak{s} \otimes_{\F_1} \F_2\).
\end{itemize}  
\end{remark}

\section{Centroid of Lie algebras}\label{sec:centroids}

Let \(L\) be a Lie algebra over the field \(\F\). The centroid is defined as
\begin{equation*}
    \Gamma(L)=\big\{\varphi \colon L \to L \mid \varphi \text{ is \(\F\)-linear and } \varphi([x,y])=[\varphi(x),y] \, \forall \, x,y \in L\big\}.
\end{equation*}
It is easy to see that \(\Gamma(L)\) is a vector space. Now we want to investigate what happens to this structure if we make a scalar extension. 
In \cite{Benkart2006}, the reader can find the following result.

\begin{theorem}[{\cite[Corollary~2.23]{Benkart2006}}]
    Let \(A,B\) be two algebras over the field \(\F\) such that \(A\) is perfect and \(B\) is unital. Then 
    \begin{equation*}
        \mathrm{Cen}(A) \otimes \mathrm{Cen}(B)\cong \big\{\chi \in \mathrm{Cen}(A \otimes B) \mid \chi \text{ has finite \(A\)-image}\big\}.
    \end{equation*}
\end{theorem}

Here we give a more general result on the dimensions of the involved algebras, where $A$ can be any type of algebra and $B$ is a finite field extension.

\begin{theorem}
    Let \(A\) be an algebra over \(\F\) and let \(\K\) be a finite extension of \(\F\). Then 
    \(\dim_\F \Gamma(A)=\dim_\K \Gamma(A \otimes_\F \K)\) where the operation in \(A \otimes_\F \K\) is defined as 
    \begin{equation*}
        (x \otimes a) \cdot (y \otimes b)=xy \otimes ab.
    \end{equation*}
\end{theorem}

\begin{proof}
 Let \(\mathcal{E}=\big\{e_1,\ldots,e_n\big\}\) be a basis for \(\K\) over \(\F\), \(\big\{v_1,\ldots,v_m\big\}\) be a basis for \(A\) over \(\F\), and \(\Gamma(A)^n\) be the Cartesian product of \(n-\)copies of \(\Gamma(A)\). The $\K$-vector space \(\Gamma(A \otimes _{\F} \K)\), in particular, is also a \(\F-\)vector space. Let us consider \(\varphi \colon \Gamma(A)^n \to \Gamma(A \otimes _{\F} \K)\), with \(\varphi((f_1,\ldots,f_n))=\tilde{f}\), where 
\begin{equation*}
    \tilde{f}(v_i \otimes 1)= \sum_{l=1}^n f_l(v_i) \otimes e_l 
\end{equation*}
and we extend by \(\K-\)linearity to every elements of \(A \otimes _{\F} \K\) (the set \(\mathcal{B}=\big\{v_i \otimes 1 \big\}_{1\leq i\leq m}\) is a \(\K-\)basis of \(L \otimes_{\F} \K\)). Firstly, we prove that $\varphi$ takes values in $\Gamma(A \otimes _{\F} \K)$.
The map $\varphi$ is well-defined, and its images $\tilde{f}$ are manifestly linear on $\K$. $\tilde{f}$ is also a centroid; indeed, we have 
\begin{align*}
    \tilde{f}(v_i \otimes 1 \cdot v_j \otimes 1)&=\tilde{f}(v_iv_j \otimes 1)\\
    &=\sum_{l=1}^n f_l(v_iv_j)\otimes e_l \\
    &=\sum_{l=1}^n f_l(v_i)v_j \otimes e_l \\
    &=\sum_{l=1}^n (f_l(v_i) \otimes e_l)\cdot(v_j \otimes 1)\\
    &=\sum_{l=1}^n (f_l(v_i) \otimes e_l)\cdot (v_j \otimes 1)\\
    &=\tilde{f}(v_i \otimes 1)\cdot(v_j \otimes 1).
\end{align*}

Similar \(\tilde{f}(v_i \otimes 1 \cdot v_j \otimes 1)=(v_i \otimes 1)\cdot\tilde{f}(v_j \otimes 1).\)The linearity of \(\varphi\) follows from the linearity of the \(f_i\)'s and \(g_j\)'s, together with the bilinearity of the tensor product. Indeed, we have

\begin{align*}
    \varphi(\alpha(f_1, \ldots,f_n)+\beta (g_1,\ldots,g_n))(v \otimes 1)&=\varphi((\alpha f_1 +\beta g_1, \ldots ,\alpha f_n +\beta g_n))(v \otimes1)\\
    &=\sum_{l=1}^n (\alpha f_l + \beta g_l)(v) \otimes e_l \\
    &= \alpha \sum_{l=1}^n f_l(v)\otimes e_l  + \beta \sum_{l=1}^n g_l(v) \otimes e_l\\
    &=\alpha \varphi((f_1,\ldots ,f_n))(v \otimes 1) \\
    &\quad+ \beta \varphi((g_1,\ldots,g_n))(v\otimes 1).
\end{align*}

Hence, \(\varphi\) is \(\F\)-linear, since its defining components are linear maps and it coincides with a linear map on a basis.

If \(\varphi((f_1,\ldots,f_n))=\tilde{f}=0\), then \(\sum_{l=1} f_l(v) \otimes e_l=0\). Then
\(f_l(v)=0\), for every \(l\) and for every \(v\) on the basis. Therefore \(f_l=0\) and this proves that $\varphi$ is injective. Lastly, we prove that \(\varphi\) is surjective. Let \(\tilde{f} \in \Gamma(A \otimes_{\F} \K)\). Hence there exist \(f_1,\ldots,f_n\) such that \(\tilde{f}(v \otimes 1)=\sum_{l=1}^n f_l(v) \otimes e_l\) for every \(v\in L\). 

Since \(\tilde{f}(x \otimes 1)+\tilde{f}(y \otimes 1)=\tilde{f}((x+y) \otimes 1)\) we have

\begin{equation*}
\sum_{l=1}^n f_l(x) \otimes e_l + \sum_{l=1}^n f_l(y) \otimes e_l = \sum_{l=1}^n f_l(x+y) \otimes e_l,
\end{equation*}
that is equivalent to
\begin{equation*}
    \sum_{l=1}^n \left(f_l(x)+f_l(y)\right) \otimes e_l = \sum_{l=1}^n f_l(x+y) \otimes e_l. 
\end{equation*}

Therefore \(f_l(x+y)=f_l(x)+f_l(y)\), for all \(1\leq l\leq n\) and for any \(x,y \in A\). Similarly, one can prove the homogeneity. 

By the identity 
\begin{equation*}
    \tilde{f}(xy \otimes 1)=\tilde{f}(x \otimes 1) \cdot (y \otimes 1)=(x \otimes 1) \cdot \tilde{f}(y \otimes 1),
\end{equation*} we obtain 
the maps \(f_l\) lie in the centroid of \(A\). Then \(\tilde{f}=\varphi((f_1,\ldots,f_n))\) because the two bilinear maps coincide for every element in $\mathcal{B}$. Hence \(\varphi\) is an \(\F-\)linear isomorphism. Now, on one hand  
\begin{equation*}
    \dim_{\F}\Gamma(A \otimes_{\F} \K)=[\K:\F]\cdot \dim_{\K}\Gamma(A \otimes_{\F} \K).
\end{equation*}
On the other hand, via the isomorphism $\varphi$, it holds
\begin{equation*}
    \dim_{\F}\Gamma(A)^n=\dim_{\F}\Gamma(A \otimes_{\F} \K).
\end{equation*}
Since $\dim_{\F}\Gamma(A)^n=n\cdot\dim_{\F}\Gamma(A)$, one obtain the statement.
\end{proof}

\begin{corollary}
    Let \(L\) be a Lie algebra over \(\F\) and let \(\K\) be a finite extension of \(\F\). Then \(\dim_\F \Gamma(L)=\dim_\K \Gamma(L \otimes_\F \K)\)
\end{corollary}

\begin{prop}
    Let \(L=s_1 \oplus \ldots s_t\) be the direct sum of perfect Lie algebras \(s_1,\ldots,s_t\) over the field \(\F\). Then \(\dim_\F \Gamma(L)=\sum_{i=1}^t \dim_\F\Gamma(s_i)\).
\end{prop}

\begin{proof}
    Let \(\varphi \in \Gamma(L)\). Then \(\varphi(s_i)=\varphi([s_i,s_i])=[\varphi(s_i),s_i]\subseteq s_i\). Hence we consider the maps \( \varphi\vert_{s_i} \colon s_i \to s_i\). We have \(\varphi\vert_{s_i} \in \Gamma(s_i)\) and \(\varphi=\sum_{i=1}^t \varphi\vert_{s_i}\).
\end{proof}

\begin{prop}
    Let \(\K\) be a Galois extension of \(\F\) of order \(n\), let \(L\) be a simple Lie algebra over \(\F\) isomorphic to a \(\K\) Lie algebra \(s\) via \(\phi\) and let \(\big\{\varphi_1,\ldots,\varphi_t\big\}\) be a \(\K\)-basis for \(\Gamma(s)\). Then a basis for \(\Gamma(L)\) is given by \(\big\{F_{ij}\big\}_{1\leq i \leq t,1 \leq j \leq p}\) defined as 
    \(F_{ij}(x)=\phi^{-1}(a_j\varphi_i(\phi(x)))\) where \(\big\{a_1,\ldots,a_n\big\}\) is basis of \(\K\) over \(\F\).
\end{prop}

\begin{proof}
    We have that \(L \otimes_\F \K \cong s \oplus \cdots \oplus s\). From the previous result \(\dim_\F \Gamma(L)=\dim_\K \Gamma(L \otimes_\F \K)=\sum_{i=1}^n \dim_\K\Gamma(s)=nt\). So we only need to verify that the set \(\big\{F_{ij}\big\}_{1\leq i \leq t,1 \leq j \leq n}\) consists of linearly independent maps.
    Let \(\sum_{i=1}^t\sum_{j=1}^n \alpha_{ij} F_{ij}(x)=0\) for all \(x\) and for some \(\alpha_{ij} \in \F\). Hence \(0=\sum_{i=1}^t\sum_{j=1}^n \alpha_{ij} \phi^{-1}(a_j\varphi_i(\phi(x)))= \phi^{-1}(\sum_{i=1}^t\sum_{j=1}^n \alpha_{ij}a_j\varphi_i(\phi(x)))\) (\(\phi^{-1}\) is \(\F\)-linear ). Since \(\phi^{-1}\) is an isomorphism this implies that \(\sum_{i=1}^t\sum_{j=1}^n \alpha_{ij}a_j\varphi_i(\phi(x)))=0.\) The map \(\phi\) is bijective and the maps \(\varphi_i\) are linearly independent over \(\K\), so we also have \(\sum_{j=1}^n \alpha_{ij} a_j=0\) for all \(i\). Lastly, the elements \(a_j\) forms a basis for \(\K\), so \(\alpha_{ij}=0\) for all \(i,j\).
\end{proof}

\section*{Acknowledgments}

Alfonso Di Bartolo and Gianmarco La Rosa were supported by the University of Palermo (FFR2026, UNIPA BsD). All authors are supported by the “National Group for Algebraic and Geometric Structures, and Their Applications” (GNSAGA — INdAM).

\section*{Conflict of Interest}
		
		The authors have no conflict of interest to declare that is relevant to this article.

\section*{Data Availability}
		Data sharing is not applicable to this article as no datasets were generated or analysed during the current study.  

\section*{Ethical Approval}

        Not applicable.   

\bibliographystyle{plain}
\bibliography{bibliography.bib}

@book{erdmann2006introduction,
	title={Introduction to Lie algebras},
	author={Erdmann, K. and Wildon, M.J.},
	series={Spriger undergraduate mathematics series},
	year={2006},
	isbn={9781846280405},
	publisher={Springer}
}

@book{humphreys2012introduction,
	title={Introduction to Lie Algebras and Representation Theory},
	author={Humphreys, J.E.},
	isbn={9781461263982},
	series={Graduate Texts in Mathematics},
	year={2012},
	publisher={Springer New York}
}

@book{jacobson1979lie,
	title={Lie algebras},
	author={Jacobson, N.},
	year={1979},
	isbn={9780486638324},
	publisher={Courier Corporation}
}

@book{knapp1996lie,
 author = {Knapp, Anthony W.},
 title = {Lie groups beyond an introduction},
 edition = {2nd ed.},
 fseries = {Progress in Mathematics},
 series = {Prog. Math.},
 issn = {0743-1643},
 volume = {140},
 isbn = {0-8176-4259-5},
 year = {2002},
 publisher = {Boston, MA: Birkh{\"a}user},
 language = {English},
 zbMATH = {1849081},
 Zbl = {1075.22501}
}

@book{Schwartz_2013, title={Quantum Field Theory and the Standard Model}, author={Schwartz, M.D.}, year={2013}, publisher={Cambridge University Press}}

@book{Gille_Szamuely_2006,
 author = {Gille, Philippe and Szamuely, Tam{\'a}s},
 title = {Central simple algebras and {Galois} cohomology},
 fseries = {Cambridge Studies in Advanced Mathematics},
 series = {Camb. Stud. Adv. Math.},
 volume = {101},
 isbn = {0-521-86103-9; 0-511-22315-3},
 year = {2006},
 publisher = {Cambridge: Cambridge University Press},
 language = {English},
 zbMATH = {5047789},
 Zbl = {1137.12001}
}

@book{Pierce,
 author = {Pierce, Richard S.},
 title = {Associative algebras},
 series = {Grad. Texts Math.},
 issn = {0072-5285},
 volume = {88},
 year = {1982},
 publisher = {Springer, Cham}
}

@article {Benkart2006,
    AUTHOR = {Benkart, Georgia and Neher, Erhard},
     TITLE = {The centroid of extended affine and root graded {L}ie
              algebras},
   JOURNAL = {J. Pure Appl. Algebra},
  FJOURNAL = {Journal of Pure and Applied Algebra},
    VOLUME = {205},
      YEAR = {2006},
    NUMBER = {1},
     PAGES = {117--145},
      ISSN = {0022-4049,1873-1376},
   MRCLASS = {17B70 (17B67)},
  MRNUMBER = {2193194},
MRREVIEWER = {Volodymyr\ Mazorchuk},
       DOI = {10.1016/j.jpaa.2005.06.007},
       URL = {https://doi.org/10.1016/j.jpaa.2005.06.007},
}

@misc{DiBartolo2026,
  author = {{Di Bartolo}, Alfonso and {Di Fatta}, Francesco Paolo and {La Rosa}, Gianmarco },
  title = {Biderivations of Lie algebras under finite field extensions},
  year = {2026},
  note = {submitted}
}

@book{Rotman1998,
 author = {Rotman, Joseph},
 title = {Galois theory.},
 edition = {2nd ed.},
 fseries = {Universitext},
 series = {Universitext},
 issn = {0172-5939},
 isbn = {0-387-98541-7},
 year = {1998},
 publisher = {New York, NY: Springer},
 language = {English},
 zbMATH = {1226184},
 Zbl = {0924.12001}
}

@book{Strade1988,
 author = {Strade, Helmut and Farnsteiner, Rolf},
 title = {Modular {Lie} algebras and their representations},
 fseries = {Pure and Applied Mathematics, Marcel Dekker},
 series = {Pure Appl. Math., Marcel Dekker},
 volume = {116},
 isbn = {0-8247-7594-5},
 year = {1988},
 publisher = {New York etc.: Marcel Dekker, Inc.},
 language = {English},
 zbMATH = {192889},
 Zbl = {0648.17003}
}

@book{Strade2004,
 author = {Strade, Helmut},
 title = {Simple {Lie} algebras over fields of positive characteristic. {I}: {Structure} theory},
 fseries = {De Gruyter Expositions in Mathematics},
 series = {De Gruyter Expo. Math.},
 issn = {0938-6572},
 volume = {38},
 isbn = {3-11-014211-2},
 year = {2004},
 publisher = {Berlin: de Gruyter},
 language = {English},
 zbMATH = {2115555},
 Zbl = {1074.17005}
}

@article{Strade2004_2,
 author = {Premet, Alexander and Strade, Helmut},
 title = {Simple {Lie} algebras of small characteristic. {IV}: {Solvable} and classical roots},
 fjournal = {Journal of Algebra},
 journal = {J. Algebra},
 issn = {0021-8693},
 volume = {278},
 number = {2},
 pages = {766--833},
 year = {2004},
 language = {English},
 doi = {10.1016/j.jalgebra.2003.10.028},
 zbMATH = {2105355},
 Zbl = {1155.17306}
}

@article{Strade2007,
 author = {Premet, Alexander and Strade, Helmut},
 title = {Simple {Lie} algebras of small characteristic. {V}: {The} non-{Melikian} case},
 fjournal = {Journal of Algebra},
 journal = {J. Algebra},
 issn = {0021-8693},
 volume = {314},
 number = {2},
 pages = {664--692},
 year = {2007},
 language = {English},
 doi = {10.1016/j.jalgebra.2007.02.059},
 zbMATH = {5201700},
 Zbl = {1139.17007}
}

@article{Strade2008,
 author = {Premet, Alexander and Strade, Helmut},
 title = {Simple {Lie} algebras of small characteristic. {VI}: {Completion} of the classification},
 fjournal = {Journal of Algebra},
 journal = {J. Algebra},
 issn = {0021-8693},
 volume = {320},
 number = {9},
 pages = {3559--3604},
 year = {2008},
 language = {English},
 doi = {10.1016/j.jalgebra.2008.08.012},
 zbMATH = {5378819},
 Zbl = {1236.17029}
}

@book{Strade2009,
 author = {Strade, Helmut},
 title = {Simple {Lie} algebras over fields of positive characteristic. {II}: {Classifying} the absolute toral rank two case},
 fseries = {De Gruyter Expositions in Mathematics},
 series = {De Gruyter Expo. Math.},
 issn = {0938-6572},
 volume = {42},
 isbn = {978-3-11-019701-3},
 year = {2009},
 publisher = {Berlin: Walter de Gruyter},
 language = {English},
 zbMATH = {5624502},
 Zbl = {1190.17001}
}

@book{Strade2013,
 author = {Strade, Helmut},
 title = {Simple {Lie} algebras over fields of positive characteristic. {III}: {Completion} of the classification},
 fseries = {De Gruyter Expositions in Mathematics},
 series = {De Gruyter Expo. Math.},
 issn = {0938-6572},
 volume = {57},
 isbn = {978-3-11-026298-8; 978-3-11-026301-5},
 year = {2013},
 publisher = {Berlin: de Gruyter},
 language = {English},
 doi = {10.1515/9783110263015},
 zbMATH = {6059163},
 Zbl = {1287.17012}
}

@book{Serre2004,
 author = {Serre, Jean-Pierre},
 title = {Corps locaux},
 edition = {4th corrected ed.},
 isbn = {2-7056-1296-3},
 year = {2004},
 publisher = {Paris: Hermann, {\'E}diteurs des Sciences et des Arts},
 language = {French},
 zbMATH = {5050432},
 Zbl = {1095.11504}
}

@book{Knus1998,
 author = {Knus, Max-Albert and Merkurjev, Alexander and Rost, Markus and Tignol, Jean-Pierre},
 title = {The book of involutions. {With} a preface by {J}. {Tits}},
 fseries = {Colloquium Publications. American Mathematical Society},
 series = {Colloq. Publ., Am. Math. Soc.},
 issn = {0065-9258},
 volume = {44},
 isbn = {0-8218-0904-0},
 year = {1998},
 publisher = {Providence, RI: American Mathematical Society},
 language = {English},
 zbMATH = {1201676},
 Zbl = {0955.16001}
}

@book{waterhouse1979,
 author = {Waterhouse, William C.},
 title = {Introduction to affine group schemes},
 fseries = {Graduate Texts in Mathematics},
 series = {Grad. Texts Math.},
 issn = {0072-5285},
 volume = {66},
 year = {1979},
 publisher = {Springer, Cham},
 language = {English},
 zbMATH = {3689557},
 Zbl = {0442.14017}
}

@article{Melville1992,
 author = {Melville, Duncan J.},
 title = {Centroids of nilpotent {Lie} algebras},
 fjournal = {Communications in Algebra},
 journal = {Commun. Algebra},
 issn = {0092-7872},
 volume = {20},
 number = {12},
 pages = {3649--3682},
 year = {1992},
 language = {English},
 doi = {10.1080/00927879208824533},
 zbMATH = {130389},
 Zbl = {0954.17502}
}

@article{block1988,
 author = {Block, Richard E. and Wilson, Robert Lee},
 title = {Classification of restricted simple {Lie} algebras},
 fjournal = {Journal of Algebra},
 journal = {J. Algebra},
 issn = {0021-8693},
 volume = {114},
 number = {1},
 pages = {115--259},
 year = {1988},
 language = {English},
 doi = {10.1016/0021-8693(88)90216-5},
 zbMATH = {4049839},
 Zbl = {0644.17008}
}

@article {Burde2021,
    AUTHOR = {Burde, Dietrich and Dekimpe, Karel and Verbeke, Bert},
     TITLE = {Almost inner derivations of {L}ie algebras {II}},
   JOURNAL = {Internat. J. Algebra Comput.},
  FJOURNAL = {International Journal of Algebra and Computation},
    VOLUME = {31},
      YEAR = {2021},
    NUMBER = {2},
     PAGES = {341--364},
      ISSN = {0218-1967,1793-6500},
   MRCLASS = {17B40},
  MRNUMBER = {4242669},
MRREVIEWER = {Zhihua\ Chang},
       DOI = {10.1142/S0218196721500181},
       URL = {https://doi.org/10.1142/S0218196721500181},
}

@article {Benkovič2024,
    AUTHOR = {Benkovič, Dominik and Eremita, Daniel},
     TITLE = {Generalized derivations of current {L}ie algebras},
   JOURNAL = {Comm. Algebra},
  FJOURNAL = {Communications in Algebra},
    VOLUME = {52},
      YEAR = {2024},
    NUMBER = {11},
     PAGES = {4603--4611},
      ISSN = {0092-7872,1532-4125},
   MRCLASS = {17B40},
  MRNUMBER = {4791644},
MRREVIEWER = {Kaiming\ Zhao},
       DOI = {10.1080/00927872.2024.2354423},
       URL = {https://doi.org/10.1080/00927872.2024.2354423},
}

@article {Zusmanovich2005,
    AUTHOR = {Zusmanovich, Pasha},
     TITLE = {Low-dimensional cohomology of current {L}ie algebras and
              analogs of the {R}iemann tensor for loop manifolds},
   JOURNAL = {Linear Algebra Appl.},
  FJOURNAL = {Linear Algebra and its Applications},
    VOLUME = {407},
      YEAR = {2005},
     PAGES = {71--104},
      ISSN = {0024-3795,1873-1856},
   MRCLASS = {17B56 (17B20 58H15)},
  MRNUMBER = {2161916},
MRREVIEWER = {Andreas\ Arvanitoyeorgos},
       DOI = {10.1016/j.laa.2005.04.018},
       URL = {https://doi.org/10.1016/j.laa.2005.04.018},
}

@article {Neeb2008,
    AUTHOR = {Neeb, Karl-Hermann and Wagemann, Friedrich},
     TITLE = {The second cohomology of current algebras of general {L}ie
              algebras},
   JOURNAL = {Canad. J. Math.},
  FJOURNAL = {Canadian Journal of Mathematics. Journal Canadien de
              Math\'ematiques},
    VOLUME = {60},
      YEAR = {2008},
    NUMBER = {4},
     PAGES = {892--922},
      ISSN = {0008-414X,1496-4279},
   MRCLASS = {17B56},
  MRNUMBER = {2432827},
MRREVIEWER = {Rolf\ Farnsteiner},
       DOI = {10.4153/CJM-2008-038-6},
       URL = {https://doi.org/10.4153/CJM-2008-038-6},
}
\end{document}